\documentclass{article}

\usepackage[english]{babel}
\usepackage{graphicx}
\usepackage{amsfonts}

\newcommand{\ebold}{{\bf e}}
\newcommand{\Mat}{{\rm Mat}}

\newcommand{\ie}{{\it i.e.\/}}
\newcommand{\R}{{\mathbb R}}
\newcommand{\N}{{\mathbb N}}

\newcommand{\X}{{\mathbb X}}
\newcommand{\Xtilde}{{\widetilde{\X}}}
\newcommand{\I}{{\mathcal I}}
\newcommand{\B}{{\mathcal B}}
\newcommand{\Btilde}{{\widetilde{\B}}}
\newcommand{\QB}{{\mathcal O}}

\newtheorem{thm}{Theorem}[section]

\newtheorem{defn}[thm]{Definition}
\newtheorem{prop}[thm]{Proposition}
\newtheorem{alg}[thm]{Algorithm}

\newenvironment{pf}{{\em Proof:\ }}{}
\newenvironment{exmp}[1][Example.]{\begin{trivlist}
\item[\hskip \labelsep {\bfseries #1}]}{\end{trivlist}}
\begin{document}

\title{Stable Border Bases for Ideals of Points}
\author{John Abbott$^*$,  Claudia Fassino\thanks{Dip. di Matematica, Universit\`a di Genova, via Dodecaneso 35, 16146 Genova, Italy - fassino@dima.unige.it}, Maria-Laura Torrente\thanks{Scuola Normale Superiore, piazza dei Cavalieri 7, 56126 Pisa, Italy - m.torrente@sns.it}}
\date{}
\maketitle

\begin{abstract}
Let~$\X$ be a set of points whose coordinates are known with limited 
accuracy; our aim is to give a characterization of the vanishing ideal~$\I(\X)$
independent of the data uncertainty.
We present a method to compute a polynomial basis~$\B$ of~$\I(\X)$
which exhibits structural stability, that is, if~$\Xtilde$ is any 
set of points differing only slightly from~$\X$, there 
exists a polynomial set~$\Btilde$ structurally similar
to~$\B$, which is a basis of the perturbed ideal~$\I(\Xtilde)$.
\end{abstract}

\section{Introduction}

In this paper we present a method for computing ``structurally stable''
border bases of ideals of points whose coordinates are affected by
errors.

If~$\X$ is a set of ``empirical'' points, representing real-world
measurements, then typically the coordinates are known only imprecisely.
Roughly speaking, if~$\Xtilde$ is another set of points, each differing by
less than the uncertainty from the corresponding element of~$\X$, then the
two sets can be considered as equivalent.  Nevertheless, it can happen that
their vanishing ideals have very different bases~---~this is a well known
phenomenon in Gr\"obner basis theory.  In order to emphasize the
``numerical equivalence'' of~$\X$ and its perturbation~$\Xtilde$, we look
for a common characterization of the vanishing ideals~$\I(\X)$
and~$\I(\Xtilde)$.  More precisely our goal is to determine a polynomial
basis~$\B$ of the vanishing ideal~$\I(\X)$ which exhibits structural
stability: namely, there is a basis~$\Btilde$ for the perturbed
ideal~$\I(\Xtilde)$, sharing the same structure as~$\B$, and whose
coefficients differ only slightly, provided that~$\Xtilde$ differs
from~$\X$ by only a small amount (up to some limit).

The decision to use border bases to describe vanishing ideals of sets of 
empirical points was due to two main reasons: border bases have always 
been considered a numerically stable tool (see~\cite{KR05}, \cite{St});
furthermore, it is easy to study their structure, \ie~the support of 
their polynomials, as it is completely determined once a suitable order 
ideal~$\QB$ has been chosen.  

We introduce the notion of stable quotient basis: given a set~$\X$ of
empirical points and a permitted tolerance~$\varepsilon$, a stable quotient
basis~$\QB$ guarantees the existence of an~$\QB$-border basis~$\Btilde$ for
the vanishing ideal~$\I(\Xtilde)$ where~$\Xtilde$ is any set of points
perturbed by amounts less than the tolerance~$\varepsilon$.  Once a stable
quotient basis~$\QB$ has been found, the corresponding stable border basis
can be obtained by some simple combinatorical and linear algebra
computations; so we focus our attention on determining~$\QB$.

An alternative approach to the problem, presented in~\cite{HKPP}, is to
use singular value decomposition of matrices to obtain a set of polynomials
which are not required to vanish on~$\X$ but must nevertheless assume
particularly small values there. In contrast, a stable border basis
always comprises polynomials which vanish on~$\X$.

\smallskip
This paper is organized as follows.  In Section~$2$ we introduce the
concepts and tools we shall use.  Section~$3$ provides a formal description
of our problem.  The main result, the SOI algorithm for computing a stable
order ideal, is presented in Section~$4$.  In Section~$5$ we give some
numerical examples illustrating the functioning of the algorithm.  Finally,
Section~$6$ is an Appendix which contains the proof of a basic result about
the first order approximation of rational functions, useful for the error
analysis of the sensitivity of the border basis computation.
 
\section{Basic definitions and notation}

This section contains basic definitions and notation used later in the
paper. To simplify the presentation, we shall implicitly suppose that each
finite set of points or polynomials is in fact a tuple, so that the 
elements are ordered in some way, and we can refer to the~$k$-th element 
using the index~$k$.

Let $n \ge 1$; we recall (see~\cite{KR00,KR05}) some basic 
concepts related to the polynomial ring $P=\R[x_1,\dots,x_n]$.
\begin{defn}\label{eval}
Let $\X= \{p_1,\dots,p_s\}$ be a non-empty finite set of points of $\R^n$ 
and let $G=\{g_1,\dots,g_k\}$ be a non-empty finite set of polynomials.
\begin{description}
\item{(a)} The ideal $\I(\X)=\{f \in P \; | \;f(p_i)=0 \; \forall p_i \in
  \X\}$ is called the {\bf vanishing ideal} of $\X$.
\item{(b)} The $\R$-linear map ${\rm{eval}}_{\X} : P \rightarrow \R^s$ 
defined by ${\rm{eval}}_{\X}(f) = (f(p_1),\dots,f(p_s))$
is called the {\bf evaluation map} associated to $\X$.
For brevity, we write $f(\X)$ to mean ${\rm eval}_{\X}(f)$.
\item{(c)} The {\bf evaluation matrix} of $G$ associated to $\X$, written
$M_{G}(\mathbb X) \in \Mat_{s,k}(\R)$, is defined as having entry $(i,j)$
equal to $g_j(p_i)$, \ie~whose columns are the images of the polynomials
$g_j$ under the evaluation map.
\end{description}
\end{defn} 

\begin{defn}
Let $\mathbb T^n$ be the monoid of power products of $P$ and let 
$\mathcal O$ be a non-empty subset of $\mathbb T^n$.
\begin{description}
\item{(a)} The {\bf factor closure} (abbr.~{\bf closure}) of $\mathcal O$ is the set $\overline{
  \mathcal O}$ of all power products in $\mathbb T^n$ which divide
  some power product of $\mathcal O$.
\item{(b)} The set $\mathcal O$ is called an {\bf order ideal} if 
$\mathcal O = \overline{\mathcal O}$,~\ie~if $\mathcal O$ is factor closed.
\item{(c)} Let $I \subseteq P$ be a zero-dimensional ideal, and $s=\dim(P/I)$; if $\mathcal O$ is factor closed and the residue classes of its elements form a basis of $P/I$ then we call it a {\bf quotient basis} for $I$.
\item{(d)} Let $\mathcal O$ be an order ideal; the {\bf border} 
$\partial \mathcal O$ of $\mathcal O$ is defined by
$$
\partial \mathcal O~=~(x_1 \mathcal O~\cup~\dots~\cup~x_n \mathcal O)~\backslash~\mathcal O
$$
\item{(e)} If $\mathcal O$ is an order ideal then the elements of the 
minimal set of generators of the monomial ideal corresponding to 
$\mathbb T^n \backslash \mathcal O$ are called the {\bf corners} of 
$\mathcal O$.
\end{description}
\end{defn}

\begin{defn}\label{borderDef}
Let $\QB=\{t_1,\dots,t_{\mu}\}$ be an order ideal, and let 
$\partial \QB = \{b_1,\dots,b_{\nu}\}$ be the border of~$\QB$.  
Let $\mathcal{B}=\{g_1,\dots,g_{\nu}\}$ be a set of polynomials having the form 
$g_j~=~b_j~-~\sum_{i=1}^{\mu} \alpha_{ij}t_i$ where each 
$\alpha_{ij} \in \R$.  Let $I \subseteq P$ be an ideal containing~$\B$.
If the residue classes of the elements of~$\QB$
form a~$\R$-vector space basis of~$P/I$ then~$\B$ is called a
{\bf border basis} of~$I$ founded on~$\QB$, or more briefly~$\B$ is an {\bf $\QB$-border basis} of~$I$.
\end{defn}

\begin{prop}{\bf (Existence and Uniqueness of Border Bases)}\\
Let~$I \subseteq P$ be a zero-dimensional ideal, and
let~$\QB=\{t_1,\dots,t_{\mu}\}$ be a quotient basis for~$I$.
Then there exists a unique $\QB$-border basis~$\B$ of~$I$.
\end{prop}

\begin{pf}
See Proposition~6.4.17 in~\cite{KR05}.
\end{pf}

\smallskip
Later on, in order to measure the distances between points of $\R^n$, we 
will use the euclidean norm $\|\cdot\|$. Additionally, given an 
$n\times n$ positive diagonal matrix $E$, we shall also use the weighted 
$2$-norm $\| \cdot \|_E$ as defined in~\cite{DBA}.
For completeness, we recall here their definitions:
\begin{eqnarray*}
\|v\|:= \sqrt{\sum_{j=1}^n v_j^2}  \qquad {\rm and } \qquad
\|v\|_E:=\|Ev\|
\end{eqnarray*}

We recall the definition of empirical point (see~\cite{St}, \cite{AFT}).
\begin{defn}\label{empirPts}
Let $p \in \R^n$ be a point and let $\varepsilon=(\varepsilon_1,\dots,
\varepsilon_n)$, with each $\varepsilon_i \in \R^+$, be the vector of the 
componentwise tolerances. An {\bf empirical point} $p^\varepsilon$
is the pair $(p,\varepsilon)$, where we call $p$ the {\bf specified value}
and $\varepsilon$ the {\bf tolerance}. 
\end{defn}

Let $p^\varepsilon$ be an empirical point.  We define its ellipsoid of
perturbations:
$$
N(p^\varepsilon) = \{\widetilde{p} \in \R^n \;:\;\|\widetilde{p}-p\|_E \le 1  \}
$$
where the positive diagonal matrix $E={\rm diag}(1/\varepsilon_1,\dots,1/
\varepsilon_n)$. This set contains all the {\bf admissible perturbations} 
of the specified value $p$, \ie~all points differing from $p$ by less than
the tolerance.

We shall assume that all the empirical points share the same 
tolerance~$\varepsilon$, as is reasonable if they derive from 
real-world data measured with the same accuracy.
In particular this assumption allows us to use the $E$-weighted norm on 
$\R^n$ to measure the distance between empirical points.

Given a finite set $\X^\varepsilon$ of empirical points all sharing the 
same tolerance $\varepsilon$, we introduce the concept of a slightly 
perturbed set of points $\widetilde{\X}$ by means of the following 
definition.
\begin{defn}
Let $\X^{\varepsilon} = \{p_1^{\varepsilon},\ldots, p_s^{\varepsilon}\}$ 
be a set of empirical points with uniform tolerance~$\varepsilon$ and 
with $\X \subset \R^n$.  Each set of points $\widetilde{\mathbb X}=
\{\widetilde{p_1},\dots,\widetilde{p_s}\} \subset \R^n$ whose elements 
satisfy
$$
(\widetilde{p_1},\dots,\widetilde{p_s}) \in \prod_{i=1}^s N(p_i^{\varepsilon})
$$
is called an {\bf admissible perturbation} of $\X^\varepsilon$.
\end{defn}

Finally we introduce the definition of distinct empirical points.
\begin{defn}
The empirical points $p_1^{\varepsilon}$ and $p_2^{\varepsilon}$, with 
specified values $p_1,p_2 \in \R^n$, are said to be {\bf distinct} if
$$
N(p_1^{\varepsilon}) \cap N(p_2^{\varepsilon}) = \emptyset
$$ 
\end{defn}

\section{The formal problem}

We shall use the concept of empirical point to describe formally the given
uncertain data: the input $\X$ is viewed as the set of specified values
of~$\X^\varepsilon$, which consists of $s$ distinct empirical points all
sharing the same fixed tolerance~$\varepsilon$.

Given the set $\X$, we want to determine a numerically stable basis~$\B$ of
the vanishing ideal~$\I(\X)$.  Intuitively, a basis~$\B$ of~$\I(\X)$ is
considered to be structurally stable if, for each admissible
perturbation~$\Xtilde$ of~$\X^\varepsilon$, it is possible to produce a
basis~$\Btilde$ of $\I(\Xtilde)$ only by means of a slight and continuous
variation of the coefficients of the polynomials of~$\B$, that is if there
exists a basis~$\Btilde$ of~$\I(\Xtilde)$ whose polynomials have the same
support as the corresponding polynomials of~$\B$.  Given a polynomial
basis~$\B$, we will call the union of the supports of its polynomials the
{\bf structure} of~$\B$.
 
A good starting point for us is the concept of border basis
(see~\cite{KR05}, \cite{St}). In fact the structure of a border basis is
easily computable and completely determined by the quotient basis~$\QB$ upon
which the border basis is founded (see Definition~\ref{borderDef}).  Using
border bases, the problem of computing a structurally stable representation
of the vanishing ideal~$\I(\X)$ thus reduces to the problem of finding a
quotient basis~$\QB$ for~$\I(\Xtilde)$ valid for every admissible
perturbation~$\Xtilde$.  The following deinition captures this notion and
generalizes it to any order ideal.
\begin{defn}
Let~$\QB$ be an order ideal, then $\QB$ is {\bf stable}
w.r.t.~$\X^\varepsilon$ if the evaluation matrix~$M_{\QB}(\Xtilde)$ has full
rank for each admissible perturbation~$\Xtilde$ of~$\X^\varepsilon$.
\end{defn} 

The following proposition highlights the importance of stable quotient bases.
\begin{prop}\label{stable}
Let~$\X^\varepsilon$ be a set of~$s$ distinct empirical points, and let~$\QB
= \{t_1,\ldots,t_s\}$ be a quotient basis for~$\I(\X)$ which is stable
w.r.t.~$\X^\varepsilon$.  Then, for each admissible perturbation~$\Xtilde$
of~$\X^\varepsilon$, the vanishing ideal~$\I(\Xtilde)$ has an~$\QB$-border
basis.  Furthermore, if $\partial{\QB}=\{b_1,\dots,b_\nu\}$ is the border
of~$\QB$ then $\Btilde$ consists of~$\nu$ polynomials of the form
\begin{eqnarray}\label{pol_border}
g_j =b_j - \sum_{i=1}^s \alpha_{ij}t_i   \qquad  {\rm for}\; j=1\dots\nu
\end{eqnarray}
where the coefficients~$a_{ij} \in \R$ satisfy
$$
b_j(\Xtilde)=\sum_{i=1}^s \alpha_{ij}t_i (\Xtilde)
$$
\end{prop}
\begin{pf}
Let~$\Xtilde$ be an admissible perturbation of~$\X^\varepsilon$ and
let~${\rm eval}_{\Xtilde}:P \rightarrow \R^s$ be the~$\R$-linear evaluation
map associated to the set~$\Xtilde$.  It is easy to prove that $\I(\Xtilde)
= \ker(\rm{eval}_{\Xtilde})$ and consequently, that the quotient ring~$P/
\I(\Xtilde)$ is isomorphic to~$\R^s$ as a vector space.  Since~$\QB$ is
stable w.r.t.~the empirical set~$\X^\varepsilon$, it follows that
$\{t_1(\Xtilde),\dots,t_s(\Xtilde)\}$ are linearly independent vectors.
Moreover~$\# \Xtilde = \#{\QB}$, so the residue classes of the elements
of~$\QB$ form a~$\R$-vector space basis of~$P/ \I(\Xtilde)$.

Let~$v_j = b_j(\Xtilde)$ be the evaluation vector associated to the power
product~$b_j$ lying in the border~$\partial \QB$; each~$v_j$ can be
expressed as
$$
v_j =\sum_{i=1}^s \alpha_{ij}t_i(\Xtilde)\qquad
{\rm for~some } \; \alpha_{ij} \in \R
$$ For each~$j$ we define the polynomial $g_j = b_j-\sum_{i=1}^s
\alpha_{ij}t_i$; by construction ${\rm eval}_{\Xtilde}(g_j)~=~0$, and so
$\Btilde~=~\{g_1,\dots,g_{\nu}\}$ is contained in~$\I(\Xtilde)$; it follows
that~$\Btilde$ is the~$\QB$-border basis of the ideal~$\I(\Xtilde)$.
\end{pf}

We observe that the coefficients $\alpha_{ij}$ of each polynomial~$g_j \in
\Btilde$ are just the components of the solution~$\alpha_j$ of the linear
system $M_{\QB}(\Xtilde)\cdot {\alpha_j} = b_j(\Xtilde)$.  It follows
that~$\alpha_{ij}$ are continuous functions of the points of the
set~$\Xtilde$ and so, since the order ideal~$\QB$ is stable
w.r.t.~$\X^\varepsilon$, they undergo only continuous variations as
$\Xtilde$ changes.  Now, the definition of stable border basis follows
naturally.

\begin{defn}
Let $\X^\varepsilon$ be a finite set of distinct empirical points, let~$\QB$
be a quotient basis for the vanishing ideal~$\I(\X)$.  If~$\QB$ is stable
w.r.t.~$\X^{\varepsilon}$ then the $\QB$-border basis~$\mathcal B$
for~$\I(\X)$ is said to be {\bf stable} w.r.t.~the set $\X^\varepsilon$.
\end{defn} 

The problem of computing a stable border basis of the vanishing ideal of a
set~$\X^\varepsilon$ of empirical points is therefore completely solved
once we have found a quotient basis~$\QB$ which is stable
w.r.t.~$\X^\varepsilon$.  If~$\QB$ exhibits these characteristics,
Proposition~\ref{stable} and the subsequent observation on the continuity
of the coefficients~$\alpha_{ij}$ prove the existence of the corresponding
stable border basis of the ideal~$\I(\X)$.  The problem of the effective
computation of a stable quotient basis is addressed in
section~\ref{stableOrder}.

\smallskip
We end this section by observing that any $\QB$-border basis of the
vanishing ideal~$\I(\X)$ is stable w.r.t.~$\X^\delta$ for a sufficiently
small value of the tolerance~$\delta$.  This is equivalent to saying that
any quotient basis~$\QB$ of~$\I(\X)$ has a ``region of stability'', as the
following proposition shows.

\begin{prop}\label{Reg_Stab}
Let~$\X$ be a finite set of points of~$\R^n$ and~$\I(\X)$ be its vanishing
ideal; let~$\QB$ be a quotient basis for~$\I(\X)$.  Then there exists a
tolerance $\delta = (\delta_1,\dots,\delta_n)$, with $\delta_i > 0$, such
that~$\QB$ is stable w.r.t.~$\X^\delta$.
\end{prop}
\begin{pf}
Let~$M_{\QB}(\X)$ be the evaluation matrix of~$\QB$ associated to the
set~$\X$; then~$M_{\QB}(\X)$ is a structured matrix whose coefficients
depend continuously on the points in~$\X$.  Since, by hypothesis,
the~$\QB$-border basis of the vanishing ideal~$\I(\X)$ exists, it follows
that~$M_{\QB}(\X)$ is invertible.  Recalling that the determinant is a
polynomial function in the matrix entries, and noting that the entries
of~$M_{\QB}(\X)$ are polynomials in the points' coordinates, we can conclude
that there exists a tolerance $\delta = (\delta_1,\dots,\delta_n)$, with
each $\delta_i > 0$, such that ${\rm det} (M_{\QB}(\Xtilde))\neq 0$ for any
perturbation~$\Xtilde$ of~$\X$.
\end{pf}

Nevertheless, since the tolerance~$\varepsilon$ of the empirical points
in~$\X^\varepsilon$ is given {\it a priori\/} by the measurements,
Proposition~\ref{Reg_Stab} does not solve our problem.  If the given
tolerance~$\varepsilon$ is larger than the ``region of stability'' of a
chosen quotient basis~$\QB$, the corresponding border basis will not be
stable w.r.t.~$\X^\varepsilon$; such a situation is shown in the following
example.
\begin{exmp}\label{line}
Let $\X^\varepsilon$ be the set of empirical points having 
$$
\X=\{(-1, -5),\;(0, -2),\;(1, 1),\;(2, 4.1) \} \subset \R^2
$$
as the set of specified values and $\varepsilon=(0.15,0.15)$ as the 
tolerance; let
$$
\widetilde{\X}=\{(-1+e_1,-5+e_2),\;(e_3,-2+e_4),\;(1+e_5,1+e_6),\;(2+e_7,4.1+e_8) \}
$$
be a generic admissible perturbation of $\X^\varepsilon$, where the 
parameters $e_i \in \R$ satisfy
$$
\|(e_1,e_2)\|_E \le 1   \qquad 
\|(e_3,e_4)\|_E \le 1   \qquad 
\|(e_5,e_6)\|_E \le 1   \qquad
\|(e_7,e_8)\|_E \le 1 
$$
Consider first $\mathcal O_1=\{1,y,x,y^2\}$, which is a quotient basis
for~$\I(\X)$.  The corresponding border basis~$\B_1$ of~$\I(\X)$ is not
stable w.r.t.~$\X^\varepsilon$.  Indeed, consider the perturbation
$\Xtilde=\{(-1,-5),\;(0,-2),\;(1,1),\;(2,4)\}$ of~$\X^\varepsilon$.  The
evaluation matrix~$M_{\QB_1}(\Xtilde)$ is singular, so no $\QB_1$-border
basis of~$\I(\Xtilde)$ exists.  It follows that~$\QB_1$ is not stable
w.r.t.~$\X^\varepsilon$ since its ``region of stability'' is too small
w.r.t.~the given tolerance~$\varepsilon$.

Now consider the quotient basis $\QB_2=\{1,y,y^2,y^3\}$, which is stable
w.r.t.~$\X^\varepsilon$.  In fact, for each perturbation~$\Xtilde$
of~$\X^\varepsilon$, we see that $M_{\QB_2}(\Xtilde)$ is a Vandermonde
matrix whose determinant is equal to
$(e_4-e_2+3)(e_6-e_2+6)(e_8-e_2+9.1)(e_6-e_4+3)(e_8-e_4+6.1)(e_8-e_6+3.1)$.
Since each~$|e_i| \le 0.15$, it follows that, for each
perturbation~$\Xtilde$, the matrix $M_{\QB_2}(\Xtilde)$ is invertible, and
so it is always possible to compute an $\QB_2$-border basis of the ideal
$\I(\Xtilde)$. \hfill$\diamondsuit$
\end{exmp}

\section{A practical solution}\label{stableOrder}

In this section we address the problem of computing an order
ideal~$\QB$ stable w.r.t.~$\X^\varepsilon$, a finite set of distinct
empirical points, and also the corresponding stable border basis when
it exists.\\

Since in real-world measurements the tolerance~$\varepsilon$ present in the
data is relatively small, our interest is focused on small
perturbations~$\Xtilde$ of the empirical set~$\X^\varepsilon$.  For this
reason our approach is based on a first order error analysis of the
problem.  We present in Section~\ref{SOI_Alg} an algorithm which computes a
stable order ideal~$\QB$ (up to first order).  In order to investigate the
stability of the order ideal~$\QB$ we use some results on the first order
approximation of rational functions (see Section~\ref{Rem_appr}) and we
introduce a parametric description of the admissible
perturbations~$\Xtilde$ of~${\X^\varepsilon}$ (see
Section~\ref{Param_Desc}).

If the output of the algorithm is actually a quotient basis then the
corresponding stable border basis~$\B$ exists for~$\I(\X)$.  To
determine~$\B$ it suffices to find the border of~$\QB$ (a simple
combinatorical computation), and then for each element of the border solve
a linear system (as described in the proof of Proposition~\ref{stable}).

\subsection{Remarks on first order approximation}\label{Rem_appr}
Let~$n \in \N$; let $F=\R(e_1,\dots,e_n)$ be the field of rational 
functions on~$\R$ and let~$f \in F$.  We use multi-index notation to 
give the Taylor expansion of~$f$ in a neighbourhood of~$0$
$$
f(e_1,\dots,e_n) =
\sum_{|\alpha| \ge 0} \frac{D^{\alpha}f(0)}{\alpha !}
e^{\alpha}
$$
We recall that given $\alpha = (\alpha_1,\dots,\alpha_n) \in \N^n$, we have
$| \alpha | = \alpha_1+\dots+\alpha_n$ and $\alpha ! = \alpha_1! \dots \alpha_n!$  Similarly $D^{\alpha} =
D_1^{\alpha_1} \dots D_n^{\alpha_n}$ (where $D_i^j = \partial^j / \partial e_i^j$) and
$e^{\alpha} = e_1^{\alpha_1} \dots e_n^{\alpha_n}$.\\
Each~$f \in F$ can be decomposed into components of homogeneous degree
in the following way:
$$
f = \sum_{k \ge 0} f_k \quad \mbox{ where }
f_k = \sum_{| \alpha| = k}\frac{D^{\alpha}f(0)}{\alpha !} e^{\alpha}
$$
and where, by convention, $D^{(0\dots 0)}f=f$.  Each polynomial 
$f_k \in \R[e_1,\dots,e_n]$ is called the {\bf homogeneous component} of~$f$ of degree~$k$.

Analogously, we can decompose a matrix $M \in \Mat_{r,s}(F)$ into
``homogeneous" parts in the following way.
\begin{defn}\label{homog}
Let $M=(m_{ij})$ be a matrix in $\Mat_{r,s}(F)$; the matrix 
$M_k~=~((m_{ij})_k)$, where~$(m_{ij})_k = \sum_{| \alpha| = k}\frac{D^{\alpha}m_{ij}(0)}{\alpha !} e^{\alpha} \in \R[e_1,\dots,e_n]$, is called
the {\bf homogeneous component} of~$M$ of degree~$k$.
\end{defn}

The following proposition characterizes the homogeneous components of 
degrees~$0$ and~$1$ of the solution and residual of a least squares
problem.
\begin{prop}\label{first_least}
Let $r,s \in \mathbb N$ with $r >s$; let~$M$ be a matrix in
$\Mat_{r,s}(F)$ and let~$v$ be a vector in $\Mat_{r,1}(F)$.  Let
$x \in \Mat_{s,1}(F)$ and  $\rho \in \Mat_{r,1}(F)$ be respectively the
solution and the residual of the least squares problem $Mx \approx v$.\\
The homogeneous components of degrees~$0$ and~$1$ of~$x$ are
\begin{eqnarray}\label{first_sol}
\begin{array}{rcl}
x_0 &=& (M_0^t M_0)^{-1} M_0^t v_0\\
x_1 &=&(M_0^t M_0)^{-1}(M_0^t v_1  + M_1^t v_0 - M_0^tM_1 x_0 - M_1^t M_0 x_0),
\end{array}
\end{eqnarray}
Moreover, the homogeneous components of degrees~$0$ and~$1$ of~$\rho$ are
\begin{eqnarray}\label{first_res}
\begin{array} {rcl}
\rho_0 &=& v_0 - M_0 x_0 \\
\rho_1 &=& v_1 - M_0 x_1 - M_1 x_0
\end{array}
\end{eqnarray}
\end{prop}
\begin{pf} This lengthy proof has been deferred to an appendix.
\end{pf}

Since the residual~$\rho$ is orthogonal to the columns of the matrix~$M$,
we have 
$$
M_0^t \rho_0 =0 \;\;\; {\rm and} \;\;\; M^t_1 \rho_0 + M^t_0 \rho_1=0
$$
but this does not imply that the vector $\rho_0+\rho_1$ is orthogonal to 
the columns of $M_0+M_1$. 
Nevertheless, since
$$ 
(M_0+M_1)^t (\rho_0+\rho_1)= M^t_0 \rho_0 + M^t_0 \rho_1 + M^t_1 \rho_0 + M^t_1 \rho_1 = M^t_1\rho_1
$$
we can assert that the vector $\rho_0 + \rho_1$ is orthogonal to the 
columns of $M_0 + M_1$, up to first order.

\subsection{A parametric description of $\X^\varepsilon$}\label{Param_Desc}

Let $\X^\varepsilon = \{p_1^\varepsilon,\ldots,p_s^\varepsilon\}$ be a
finite set of distinct empirical points with specified values~$\X \subset
\R^n$; we represent an admissible perturbation
$\Xtilde=\{\widetilde{p}_1,\ldots,\widetilde{p}_s\}$ of
${\X^\varepsilon}$ by using first order infinitesimals for the 
perturbation in each coordinate. In particular we express~$\Xtilde$
as a function of~$ns$ variables
$$
\ebold = (e_{11},\dots,e_{s1},e_{12},\dots,e_{s2},\dots,e_{1n},\dots,e_{sn})
$$
called {\bf error variables}; specifically, we have
\begin{eqnarray*}
\widetilde{p_k} =
\left (
p_{k1} + e_{k1},\; p_{k2} + e_{k2}, \dots p_{kn}  +  e_{kn} \right)
\end{eqnarray*}
The condition that each $\widetilde{p}_k$ is an admissible perturbation of
the point~$p_k$ is equivalent to the following:
\begin{eqnarray}\label{constraint}
\|(e_{k1},\dots,e_{kn})\|_E \le 1
\end{eqnarray}

We denote by
$\Xtilde(\ebold)=(\widetilde{p}_1(\ebold),\dots,\widetilde{p}_s(\ebold))$
a generic admissible perturbation of~$\X^\varepsilon$.  We observe that
the coordinates of each perturbed point~$\widetilde{p}_k(\ebold)$ are 
elements of the polynomial ring~$R=\R[\ebold]$ and that each variable~$e_{kj}$ represents the perturbation in the~$j$-th coordinate of the 
original point~$p_k$. 
The domain of the perturbed set ${\Xtilde}(\ebold)$, viewed as a 
function of~$ns$ variables, is denoted by~$D_\varepsilon$. 
Obviously, if $\ebold \in D_\varepsilon$ we have
$$ 
\|\ebold\|^2 =\sum_{j=1}^n\sum_{k=1}^s e_{kj}^2 \le \sum_{j=1}^ns \varepsilon_j^2,
$$
and consequently
\begin{eqnarray}\label{dominio}
\| \ebold \| \le \sqrt{s} \|\varepsilon \|
\end{eqnarray}

To keep evident the dependence on the error variables~$\ebold$, we extend
the concepts of evaluation map of a polynomial~$f \in P$ and evaluation
matrix of a set of polynomials $G=\{g_1,\dots,g_k\} \subset P$ (see
Definition~\ref{eval}) to a generic perturbed set 
$\Xtilde(\ebold)$, using the following notation:
\begin{eqnarray*}
  {\mbox{eval}}_{\Xtilde}(\ebold)(f)= \left( f(\widetilde p_1(\ebold)),\dots,f(\widetilde p_s(\ebold)) \right)  \in R \times \dots \times R = R^s
\end{eqnarray*}
for brevity denoted by~$f(\Xtilde(\ebold))$; similarly we write the
evaluation matrix
\begin{eqnarray*}
M_G(\Xtilde(\ebold)) =
\left( g_1(\widetilde{\X}(\ebold)),\dots, g_k(\widetilde{\X}(\ebold)) \right)
\end{eqnarray*}

\subsection{The SOI Algorithm}\label{SOI_Alg}

In this section we present the SOI algorithm which computes, up to
first order, an order ideal~$\QB$ stable w.r.t~the empirical
set~$\X^\varepsilon$.

The strategy for computing a stable order ideal~$\QB$ is the following.  As
in the Buchberger-M\"oller algorithm~(\cite{BM}, \cite{Ab00}) the order
ideal~$\QB$ is built stepwise: initially~$\QB$ comprises just the power
product~$1$; then at each iteration, a new power product~$t$ is considered.
If the evaluation matrix $M_{\QB \cup \{t\}} (\Xtilde(\ebold))$ has full
rank for all~$\ebold$ in~$D_{\varepsilon}$ then~$t$ is added to~$\QB$;
otherwise~$t$ is added to the corner set of the order ideal.

A first observation concerns the choice of the power product~$t$ to analyze
at each iteration: any strategy that chooses a term~$t$ such that the
set~$\QB \cup \{t\}$ preserves the property of being an order ideal can be
applied.  A possible technique is the one used in the Buchberger-M\"oller
algorithm, where the power product~$t$ is chosen according to a fixed term
ordering~$\sigma$. The version of the SOI Algorithm presented below employs
this latter stategy. Note that~$\sigma$ is only used as a computational
tool for choosing~$t$; in fact the final computed set~$\mathcal O$ is not,
in general, the same as that which would be obtained processing the
set~$\X$ by the Buchberger-M\"oller algorithm with the same term ordering.

Another observation concerns the main check of the algorithm: note that the
rank condition is equivalent to checking whether~$\rho(\ebold)$, the
component of the evaluation vector~$t(\Xtilde(\ebold))$ orthogonal to the
column space of the matrix $M_{\QB}(\Xtilde(\ebold))$, vanishes for
any~$\ebold \in D_\varepsilon$.  This check is greatly simplified by our
restriction to first order error terms.

\begin{alg}\label{SOI}
\emph{(Stable Order Ideal Algorithm)}\\
Let~$\sigma$ be a term ordering on~$\mathbb T^n$ and let
$\X^\varepsilon=\{p_1^\varepsilon,\dots,p_s^\varepsilon\} $ be a finite set
of distinct empirical points, with~$\X \subset \R^n$ and a common tolerance
$\varepsilon=(\varepsilon_1,\dots,\varepsilon_n)$.  Let $\ebold
=(e_{11},\dots,e_{sn})$ be the error variables whose constraints are given
in~$(\ref{constraint})$.  Consider the following sequence of instructions.
\begin{description}
\item[S1] Start with the  lists $\QB=[1]$, $L=[x_1,\dots,x_n]$, 
the empty list  $C=[\;]$, and the matrices
$M_0 \in {\rm Mat}_{s,1}(\R)$ with all the elements equal to $1$, and 
$M_1 \in {\rm Mat}_{s,1}(R)$ with all the elements equal to $0$.

\item[S2] If $L=[\;]$ then return the set~$\QB$ and stop. 
Otherwise let $t =\min_{\sigma}(L)$ and delete it from~$L$.

\item[S3] Let~$v_0$ and~$v_1$ be the homogeneous components of degrees~$0$ and~$1$ of the evaluation vector $v = t(\Xtilde(\ebold))$.
Solve up to first order the least squares problem
$M_{\QB} (\Xtilde(\ebold)) \; \alpha(\ebold) \approx v$, by computing
the vectors
\begin{eqnarray*}
\rho_0 &=& v_0 - M_0 \alpha_0\\
\rho_1 &=& v_1 - M_0 \alpha_1 - M_1 \alpha_0
\end{eqnarray*}
where
\begin{eqnarray*}
\alpha_0 &=& (M_0^t M_0)^{-1} M_0^t v_0\\
\alpha_1 &=&(M_0^t M_0)^{-1}(M_0^t v_1  + M_1^t v_0 - M_0^t M_1 \alpha_0 - M_1^t M_0 \alpha_0).
\end{eqnarray*}

\item[S4] Let $C_t \in {\rm Mat}_{s,sn}(\R)$ be such that
$\rho_1=C_t \ebold $. Compute the minimal 2-norm solution~$\hat \ebold$ 
of the underdetermined system $C_t \ebold= -\rho_0$ \cite{DH93}.

\item[S5] If $ \| \hat \ebold\| > \sqrt{s}\| \varepsilon \|$ then adjoin 
the vector~$v_0$ as a new column of~$M_0$ and the vector~$v_1$ as a new 
column of~$M_1$.  
Append the power product~$t$ to~$\QB$, and add to~$L$ those 
elements of $\{x_1 t,\dots,x_n t\}$ which are not multiples of an 
element of~$L$ or~$C$.  Continue with step S2.

\item[S6] Otherwise append~$t$ to the list~$C$, and remove from~$L$ all
multiples of~$t$.  Continue with step S2.
\end{description}
\end{alg}

\begin{thm}
Algorithm \ref{SOI} stops after finitely many steps and returns a set
$\QB \subset \mathbb{T}^n$ which is an order ideal stable (up to 
first order) w.r.t.~the empirical set~$\X^\varepsilon$.
Furthermore, if~$\# \QB = s$ then~$\I(\X)$ has a
corresponding stable border basis w.r.t.~$\X^\varepsilon$.
\end{thm}

\begin{pf}
First we claim that the vectors~$\rho_0$, $\rho_1$, $\alpha_0$, $\alpha_1$
computed in step~$S3$ are the homogeneous components of degrees~$0$ and~$1$
of the residual~$\rho(\ebold)$ and of the solution~$\alpha(\ebold)$ to the
least squares problem
\begin{eqnarray}\label{LSP}
M_{\QB}(\Xtilde(\ebold)) \; \alpha(\ebold) \approx
t(\Xtilde(\ebold))
\end{eqnarray}
where~$t$ is the power product being considered at the current iteration,
and~$\QB$ is the order ideal computed so far.  To prove this claim it is
sufficient to apply Proposition~\ref{first_least} to~$(\ref{LSP})$ and to
observe that the matrices~$M_0$ and~$M_1$ coincide with the homogeneous
components of degrees~$0$ and~$1$ of~$M_{\QB}(\Xtilde(\ebold))$.  Clearly,
this is true at the first iteration, since $M_{\QB}(\Xtilde(\ebold))=(1
\dots 1)^t$.  We apply induction on the number of iterations.  Assume
that~$M_0$ and~$M_1$ are the components of degrees~$0$ and~$1$
of~$M_{\QB}(\Xtilde(\ebold))$ and suppose that the power product~$r$ is
added to~$\QB$.  Since the last column of $M_{\QB \cup
  \{r\}}(\Xtilde(\ebold))$ is given by~$r(\Xtilde(\ebold))$, whose
components of degrees~$0$ and~$1$ are~$r_0$ and~$r_1$ respectively, the new
matrices~$[M_0,\; r_0]$ and~$[M_1,\; r_1]$ are the components of
degrees~$0$ and~$1$ of $M_{\QB \cup \{r\}}(\Xtilde(\ebold))$.  We conclude
that the vectors $\rho_0+\rho_1(\ebold)$ and $\alpha_0+\alpha_1(\ebold)$
coincide with~$\rho(\ebold)$ and~$\alpha(\ebold)$, up to first order.

\smallskip
Now we prove the finiteness and the correctness of Algorithm~\ref{SOI}.\\
First we show finiteness.  At each iteration the algorithm performs either
step~S5 or step~S6.  We observe that step~S5 can be executed at most~$s-1$
times; in fact, when~$M_0$ becomes a square matrix,~{\it i.e.} after~$s-1$
iterations of step~S5, the residual vector~$\rho_0$ is zero, and
consequently the minimal $2$-norm solution~$\hat{\ebold}$ of the linear
system $C_t \ebold~=~-\rho_0$ is also zero.  Moreover, step~S5 is the only
place where the set~$L$ is enlarged with a finite number of terms, while
each iteration removes from~$L$ at least one element; we conclude that the
algorithm reaches the condition~$L = [\;]$ after finitely many iterations.

\smallskip
In order to show correctness we prove, by induction on the number of
iterations and using a first order error analysis, that the output
set~$\QB$ is an order ideal stable w.r.t.~$\X^\varepsilon$.  This is
clearly true after zero iterations, \ie~after step~S1 has been executed.
By induction assume that a number of iterations has already been performed
and that the set~$\QB$ satisfies the given requirements; let us follow the
steps of the new iteration, in which a power product~$t$ is considered.  If
step~S6 is performed the claim is true because~$\QB$ does not change.
Otherwise, if step~S5 is performed, the set~$\QB^*= \QB \cup \{t\}$ is an
order ideal by construction.  Further, since the minimal $2$-norm
solution~$\hat \ebold$ of the linear system $C_t \ebold = - \rho_0$
satisfies condition $\|\hat \ebold\| > \sqrt{s} \; \|\varepsilon\|$ it
follows that~$\hat \ebold$ does not belong to~$D_\varepsilon$ and that the
vector $\rho_0+\rho_1(\ebold)$ does not vanish as~$\ebold$ varies
in~$D_\varepsilon$.  Therefore, up to first order, we can
consider~$\rho(\ebold)$ as a non-vanishing vector for each
perturbation~$\Xtilde(\ebold)$, \ie~the matrix~$M_{\QB^*}(\Xtilde(\ebold))$
has full rank for each~$\ebold \in D_\varepsilon$.

\smallskip
For the last part of the theorem we simply observe that when
$\#\QB = s$ then $\QB$ is a quotient basis; the rest is immediate.
\end{pf}

\section{Numerical examples}

In this section we present some numerical examples to show the 
effectiveness of the SOI algorithm. 
Our algorithm is implemented using the C++ language and the 
CoCoALib, see \cite{Co}, and all computations have been performed on an Intel 
Pentium~M735 processor (at~1.7~GHz) running GNU/Linux. 
In all the examples, the SOI algorithm is performed using a fixed 
precision of $1024$ bits for the RingTwinFloat implemented in CoCoALib, 
and the degree lexicographic term ordering~$\sigma$; in addition, the 
coefficients of the polynomials are displayed as truncated decimals. 

The first two examples show how the SOI algorithm detects the simplest 
geometrical configuration almost satisfied by the empirical 
set~$\X^\varepsilon$.

\begin{exmp}{\bf Almost aligned points}\\
We consider the empirical set~$\X^\varepsilon$ given in 
Example~\ref{line}; we recall here the points in~$\X$
$$
\X=\{(-1, -5),\;(0, -2),\;(1, 1),\;(2, 4.1) \} \subset \R^2
$$  
and the tolerance $\varepsilon = (0.15, 0.15)$.\\
Applying algorithm SOI to~$\X^\varepsilon$ we obtain the quotient basis
$\QB = \{1,y,y^2,y^3\}$ which is
stable w.r.t.~$\X^\varepsilon$, as we proved in Example~\ref{line}.
As~$\QB$ is a quotient basis we can compute the border basis
founded on it:
\begin{eqnarray*}
\mathcal B =
\left \{
\begin{array}{rcl}
x &+& 0.0002y^3 + 0.0012 y^2 - 0.3328y -0.6686 \\ 
xy &+& 0.0008y^3 - 0.3286y^2 - 0.6643y - 0.0079\\ 
xy^2 &-& 0.3301y^3 - 0.6471y^2 + 0.0098y - 0.0326\\ 
xy^3 &-& 0.0199y^3 - 7.1199y^2 - 7.3933y + 13.533\\
y^4 &+& 1.9 y^3 - 21.6y^2 - 22.3y + 41
\end{array}
\right.
\end{eqnarray*}
Note that the lowest degree polynomial of~$\B$, 
$f = x+ 0.0002y^3 + 0.0012 y^2 - 0.3328y -0.6686$, highlights the fact 
that~$\X$ contains ``almost aligned'' points.
In fact, if we neglect the terms with smallest coefficients,~$f$ simplifies
to $x - 0.3328y -0.6686$.  Since the coefficients of
a polynomial are continuous functions of its zeros and the quotient basis~$\QB$ is stable w.r.t.~$\X^\varepsilon$, we can conclude that there 
exists a small perturbation~$\Xtilde$ of~$\X$ containing aligned 
points and for which the associated evaluation matrix~$M_{\QB}(\Xtilde)$ is invertible. 
A simple example of such a set is given by 
$\Xtilde = \{(-1,-5),\;(0, -2),\;(1, 1),\;(2, 4)\}$. 

\smallskip
A completely different result is obtained by applying to the set~$\X$ the 
Buchberger-M\"oller algorithm w.r.t.~the same term ordering~$\sigma$.
The $\sigma$-Gr\"obner basis~$\mathcal G$ of~$\I(\X)$~is:
\begin{eqnarray*}
\mathcal G =
\left \{
\begin{array}{rcl}
x^2 &-& 1/9y^2 - 121/30x + 9/10y + 101/45\\
xy &-& 1/3y^2 - 41/10x + 7/10y + 41/15\\
y^3 &+& 6y^2 + 516243/100x - 171781/100y - 172581/50
\end{array}
\right.
\end{eqnarray*}
and the associated quotient basis is $\QB_{\sigma}(\I(\X)) = \mathbb{T}^2
\backslash \rm{LT}_{\sigma} \{\I(\X)\} = \{1,y,x,y^2\}$.  We observe
that~$\QB_{\sigma}(\I(\X))$ is not stable (see Example~\ref{line}) because
the evaluation matrix~$M_{\QB_{\sigma}}(\Xtilde)$ is singular for some
admissible perturbations of~$\X$.  In particular, the information that the
points of~$\X$ are ``almost aligned'' is not at all evident from~$\mathcal
G$.
\end{exmp}

\begin{exmp}{\bf Empirical points close to an ellipse}\\
Let $\X \subset \R^2$ be a set of points created by perturbing by less than 
$0.1$ the coordinates of $10$ points lying on the ellipse 
$x^2+0.25 y^2 -25=0$,
\begin{eqnarray*}
\X&=&\{( -5.07, 0.02),(4.98,0),(3.05, 8.07), (3.01,-8.02),(-3.02,7.99),\\
&\;&(-2.98,-8),(4.01,5.94), (3.98,-6.06),(-3.92,6.03), (-4.01,-6)\}
\end{eqnarray*}
Let~$\X^\varepsilon$ be the set of empirical points whose set of specified
values is~$\X$ and whose common tolerance is $\varepsilon = (0.1, 0.1)$.\\
Applying SOI on~$\X^\varepsilon$ we obtain, after~$11$ iterations, the 
stable quotient basis
$$ 
\QB = \{1, y, x, y^2, xy, y^3, xy^2, y^4, xy^3, xy^4\}
$$
We use linear algebra to compute the corresponding stable 
border basis~$\B$ of~$\I(\X)$.  We can
identify the ``almost elliptic'' configuration of the points of~$\X$
by looking at~$f$ the lowest degree polynomial contained in~$\B$:
\begin{eqnarray*}
f &=& x^2 + 0.273 y^2 - 25.250 + 10^{-2}
(0.004xy^4 + 0.020xy^3 - 0.034y^4 - 0.489xy^2 \\
&\;& -0.177y^3 - 1.371xy + 9.035x + 9.810y)
\end{eqnarray*}
We observe that~$f$ highlights the fact that~$\X$ contains points close to
an ellipse.  In fact, if we neglect the terms with smallest
coefficients,~$f$ simplifies to $x^2 + 0.273 y^2 - 25.250$.  Since the
coefficients of a polynomial are continuous functions of its zeros and the
quotient basis~$\QB$ is stable w.r.t.~$\X^\varepsilon$, we can conclude
that there exists a small perturbation~$\Xtilde$ of~$\X$ containing points
lying on an ellipse and such that the associated evaluation
matrix~$M_{\QB}(\Xtilde)$ is invertible.  A simple example of such a set is
given by
\begin{eqnarray*}
\widetilde{\X} &=&
\{( -5, 0),(5,0),(3, 8), (3,-8),(-3,8),\\
& &(-3,-8), (4,6),(4,-6),(-4,6),(-4,-6)\}
\end{eqnarray*}
\end{exmp}

\smallskip
\begin{exmp}{\bf Empirical points close to a circle}\\
In this example we show the behaviour of the SOI algorithm when 
applied to several sets of points with similar geometrical 
configuration but with different cardinality.\\
Let $\X_1, \X_2, \X_3, \X_4 \subset \R^2$ be sets of points created by
perturbing by less than $0.01$ the coordinates of $8, 16, 32$ and $64$
points lying on the circumference $x^2+y^2 -1=0$, and let $\varepsilon =
(0.01, 0.01)$ be the tolerance.  We summarize in Table~\ref{tab1} the
numerical tests performed by applying the SOI algorithm to the empirical
set~$\X_i^\varepsilon$, for $i=1 \dots 4$.  The first two columns of the
table contain the name of the processed set and the value of its
cardinality.  The column labelled with ``Corners'' refers to the set of
corners of the stable order ideal computed by the algorithm; the column
labelled with ``Time'' contains the time taken to compute the quotient bases.

\begin{table}[h]
\begin{center}
\begin{tabular}{|c|c|c|c|}
\hline
Input & $\# \X_i$ & Corners & Time\\
\hline
$\X_1$ & 8 & $\{x^2, xy^3, y^5\}$ & 0.5 s\\
$\X_2$ & 16 & $\{x^2, xy^7, y^9\}$ & 8.5 s\\
$\X_3$ & 32 & $\{x^2, xy^{15}, y^{17}\}$ & 79 s\\
$\X_4$ & 64 & $\{x^2, xy^{31}, y^{33}\}$ & 2320 s\\
\hline
\end{tabular}
\caption{SOI on sets of points close to a circle}\label{tab1}
\end{center}
\end{table}
Note that the set of corners of the stable quotient bases computed by the 
SOI algorithm always contain the power product~$x^2$: this means that there is an
``almost linear dependence'' among the power products  
$\{1, y, x, y^2, xy, x^2\}$ and that some useful information on 
the geometrical configuration of the points could be found.
\end{exmp}

\begin{exmp}{\bf Empirical points close to an hyperbola, a circle and 
a cubic}\\

The numerical tests suggest that in most cases the SOI algorithm 
computes a stable quotient basis, allowing us to determine a stable border basis of~$\I(\X)$.  Nevertheless, this is not true in general, as the following
example illustrates.
Let~$\X^\varepsilon$ be the set of distinct empirical points having 
$$
\X = \{(1,6), (2,3), (2.449,2.449), (3,2), (6,1)\} \subset \R^2
$$ 
as the set of specified values and $\varepsilon =(0.25, 0.25)$ as the 
tolerance.\\
Applying the algorithm SOI to~$\X^\varepsilon$, we obtain the stable 
order ideal~$\QB =\{1,y,x,y^2\}$; however, this is not a quotient
basis, so we cannot obtain a corresponding stable border basis.
This is due to the fact that the points of~$\X$ lie close to the hyperbola
$xy - 6 = 0$, the circle $(x-6)^2+(y-6)^2-25 = 0$ and the cubic
$y^3-12y^2+6x+47y-73 = 0$.  So, if the tolerance~$\varepsilon$ is too big,
they ``almost satisfy'' all of them.

Observe how the problem does not arise if we use a smaller tolerance,
{\it e.g.} $\delta = (0.2, 0.2)$.
Applying SOI to~$\X^{\delta}$ we obtain the stable quotient basis
$\QB'=\{1,y,x,y^2,y^3\}$, and its corresponding border basis:
\begin{eqnarray*}
\mathcal B' = \left\{
\begin{array}{rcl}
xy  &+& 0.0047y^3 - 0.0560y^2 + 0.0280x + 0.2194y - 6.336\\ 
x^2 &-& 0.4265y^3 + 6.118y^2 - 14.559x - 32.047y + 77.711\\
xy^2 &+& 0.0114y^3 - 0.1372y^2 + 0.0686x - 5.463y - 0.8231\\
y^4 &-& 14.477y^3 + 76.724y^2 - 14.862x - 188.419y + 214.345\\
xy^3 &+& 0.0280y^3 - 6.336y^2 + 0.1680x + 1.316y - 2.016
\end{array}
\right.
\end{eqnarray*}
\end{exmp}

\section{Appendix}\label{Appendix}

In this section we present the proof of Proposition~\ref{first_least}.
\begin{pf}
Let~$n \in \N$; first we prove a result on the homogeneous components of 
degrees~$0$ and~$1$ of the inverse of a square matrix 
$A \in \Mat_{n,n}(F)$.\\
Let~$A$ be a non singular matrix in~$\Mat_{n,n}(F)$ and let~$B$ be the 
inverse of~$A$. The homogeneous components~$B_0$ and~$B_1$ of~$B$ are
given by
\begin{eqnarray}\label{inversa}
B_0 = A_0^{-1}  \qquad \qquad B_1 = - A_0^{-1} A_1  A_0^{-1}
\end{eqnarray}
Define $\Delta A = \sum_{i \ge 2} A_i$ and 
$\Delta B = \sum_{i \ge 2} B_i$, so we have
$$ 
A=A_0+A_1+\Delta A \quad {\rm and } \quad B=B_0+B_1+\Delta B
$$
Since $A B= I$, where~$I$ is the~$n \times n$ identity matrix, we have
$$ (A_0+A_1+\Delta A)(B_0+B_1+\Delta B)= I$$
and our claim is immediate.\\
Now we prove the result of the proposition.  Applying the classical least 
squares method to the linear system $Mx \approx v$ we obtain
\begin{eqnarray}
x &=& (M^t M)^{-1} M^t v \label{soluz}\\
\rho &=& v - Mx \label{residuo}
\end{eqnarray}
Applying to (\ref{residuo}) the homogeneous degree decomposition up to degree~$1$ we have
$$
\rho_0 + \rho_1 = (v_0 - M_0 x_0) + (v_1 - M_0 x_1 - M_1 x_0)
$$
thus~(\ref{first_res}) follows.\\
Since $(M^t M)_0=M_0^t M_0$ and $(M^t M)_1=M_0^tM_1 +M_1^t M_0$, from formula~(\ref{inversa}) we have to first order,
$$
(M^tM)^{-1} \cong (M_0^t M_0)^{-1}- (M_0^t M_0)^{-1}( M_0^tM_1 + M_1^tM_0)(M_0^t M_0)^{-1}
$$
Up to degree~$1$, formula~$(\ref{soluz})$ becomes
\begin{eqnarray*}
x_0&+&x_1 = (M^t M)^{-1}_0 (M_0^t v_0 + M_0^t v_1  + M_1^t v_0) - (M^tM)^{-1}_1 M_0^t v_0\\
&=& (M_0^t M_0)^{-1}\Big(M_0^t v_0 + M_0^t v_1  + M_1^t v_0 - (M_0^tM_1 + M_1^t M_0)
(M_0^tM_0)^{-1} M_0^tv_0 \Big )
\end{eqnarray*}
and so
\begin{eqnarray*}
x_0 &=&  (M_0^t M_0)^{-1}M_0^t v_0\\
x_1 &=& (M_0^t M_0)^{-1}(M_0^t v_1  + M_1^t v_0 - M_0^tM_1 x_0 - M_1^t M_0 x_0)
\end{eqnarray*}
thus the proof is concluded.
\end{pf}

\bigskip
{\bf Acknowledgements.}
The authors would like to thank Prof.~L.~Robbiano for his useful and 
constructive remarks.

During the development of this work John Abbott was a member of a
project financially supported by the Shell Research Foundation.

This work was partly conducted during the Special Semester on
Gr\"obner Bases (from 1st February to 31st July 2006) organized by 
RICAM (Radon Institute for Computational and Applied Mathematics) of the
Austrian Academy of Sciences and RISC (Research Institute for Symbolic 
Computation) of the Johannes Kepler University, Linz, Austria, under the
direction of Professor Bruno Buchberger.


\begin{thebibliography}{}
\bibitem{Ab00}
Abbott, J., Bigatti, A., Kreuzer, M., Robbiano, L. (2000).
Computing ideals of points. \emph{J.~Symb. Comput.}, 30:341-356.

\bibitem{AFT}
Abbott, J., Fassino, C., Torrente, M. (2007).
Thinning Out Redundant Empirical Data.
\emph{Mathematics  in Computer Science} - To appear.

\bibitem{BM}
Buchberger, B., M\"oller, H.~M. (1982).
The construction of multivariate polynomials with preassigned zeros.
\emph{Proc. EUROCAM '82, Lecture Notes in Comp.Sci.}, 144:24-31.

\bibitem{Co}
CoCoA Team.
CoCoA: a system for doing computations in Commutative Algebra.
Available at \verb|http://cocoa.dima.unige.it/|

\bibitem{DBA}
Dahlquist, G., Bj\"orck, \AA., Anderson, N. (1974).
\emph{Numerical Methods}. Englewood Cliffs, New Jersey.

\bibitem{DH93}
Demmel, J.~W., Higham, N.~J. (1993).
Improved error bounds for underdetermined system solvers.
\emph{SIAM J.~Matrix Anal. Appl.}, 14:1-14.

\bibitem{HKPP}
Heldt, D., Kreuzer, M., Pokutta, S., Poulisse, H. (2006).
Approximate computation of zero-dimensional polynomial ideals.
Available at \verb|http://www.mathematik.uni-dortmund.de/algebraic-oil/|.

\bibitem{KR00}
Kreuzer, M., Robbiano, L. (2000).
\emph{Computational Commutative Algebra 1}. Springer, Berlin.

\bibitem{KR05}
Kreuzer, M., Robbiano, L. (2005).
\emph{Computational Commutative Algebra 2}. Springer, Berlin.

\bibitem{Sauer}
Sauer, T. (2007).
Approximate varieties, approximate ideals and dimension reductions.
\emph{Numerical Algorithms}. To appear.

\bibitem{St}
Stetter, H. (2004)
\emph{Numerical Polynomial Algebra}. SIAM, Philadelphia.
\end{thebibliography}
\end{document}